\documentclass[10pt,twoside,leqno]{my_cras_i-fr}
\usepackage{amssymb,amsbsy,amsmath,amsfonts,amssymb,amscd}
\usepackage[english,francais]{babel}
\usepackage{times}
\newcommand{\C}                 {\mathbb{C}}

\newtheorem{e-proposition}{Proposition\it}

\newtheorem{e-definition}{Definition\rm}

\newtheorem{theoreme}{Th\'eor\`eme\it}

\newtheorem{proposition}{Proposition\it}
\newtheorem{corollaire}{Corollaire\it}

\renewcommand{\theequation}{\arabic{equation}}
\ISSN{????}
\PIT{}
\PIT{FLA}
\PXMA{????}
\Add{?}

\Volume{}
\Year{2002}
\AuteurCourant{H.~Gaussier et J.~Merker}
\TitreCourant{Alg\'ebrisabilit\'e de sous-vari\'et\'es
analytiques r\'eelles de $\C^n$}
\Journal
\Rubrique{G\'eom\'etrie diff\'erentielle}{Differential Geometry}
\SousRubrique{}{}
\PresentePar{\'Etienne}{Ghys}
\Recu{ 11 octobre 2002}{}{}

\begin{document}\label{firstpage}


\title{Sur l'alg\'ebrisabilit\'e locale de sous-vari\'et\'es\\
analytiques r\'eelles g\'en\'eriques de $\C^n$
}

\author{%
Herv\'e GAUSSIER~$^{\text{a}}$,\ \
Jo\"el MERKER~$^{\text{a}}$.\ \
}

\address{%
\begin{itemize}\labelsep=2mm\leftskip=-5mm
\item[$^{\text{a}}$]
LATP, UMR 6632, Universit\'e de Provence, 39 rue Joliot-Curie, 13453
Marseille cedex 13, France
\item[]
Courriel~:
\{gaussier,merker\}@cmi.univ-mrs.fr
\end{itemize}
}

\maketitle

\thispagestyle{empty}

\begin{Resume}{On \'etablit que le (pseudo)groupe local des
biholomorphismes stabilisant une sous-vari\'et\'e alg\'ebrique
r\'eelle, minimale, finiment non d\'eg\'en\'er\'ee de $\mathbb C^n$,
est un groupe de Lie local alg\'ebrique r\'eel. On en d\'eduit des
conditions n\'ecessaires pour l'alg\'ebrisabilit\'e locale de tubes
analytiques r\'eels rigides de codimension quelconque dans $\mathbb
C^n$.
}\end{Resume}

\begin{Etitle}{%
On the local algebraizibility of real analytic generic
submanifolds in $\mathbb C^n$}
\end{Etitle}

\begin{Abstract}{%
We prove that the local (pseudo)group of biholomorphisms
stabilizing a minimal, finitely nondegenerate real algebraic
submanifold in $\mathbb C^n$ is a real algebraic local Lie group. We
deduce necessary conditions for the local algebraizability of real
analytic rigid tubes of arbitrary codimension in $\mathbb C^n$.}
\end{Abstract}

\setcounter{section}{0}
\section{Introduction}

Une sous-vari\'et\'e analytique r\'eelle de $\mathbb C^n$ est dite
{\sl alg\'ebrique} (au sens de Nash, {\sl cf}.~\cite{BCR98}) si elle
est d\'efinie localement par l'annulation de polyn\^omes r\'eels et
{\sl localement alg\'ebrisable} en un de ses points $p$ s'il existe un
syst\`eme de coordonn\'ees holomorphes locales centr\'e en $p$ dans
lequel elle est alg\'ebrique. S'il est clair que toute
sous-vari\'et\'e analytique r\'eelle de dimension $k$ dans $\mathbb
R^{2n}$ est localement \'equivalente, par une transformation
analytique r\'eelle, \`a un plan de dimension $k$,
il n'en va pas de m\^eme pour l'alg\'ebrisabilit\'e au moyen d'une transformation holomorphe. X.~Huang, S.~Ji et S.T.~Yau ont r\'ecemment
pr\'esent\'e dans \cite{HJY01} le premier exemple explicite
$\{(z,w) \in \mathbb C^2 : {\rm Im}\, w = \exp(z\bar{z})\}$ d'une
hypersurface analytique r\'eelle Levi non d\'eg\'en\'er\'ee 
de $\mathbb C^2$, alg\'ebrisable en
aucun de ses points. La d\'emonstration dans \cite{HJY01} repose sur
la m\'ethode dite d'\'equivalence due \`a \'E.~Cartan \cite{Ca32} et
d\'evelopp\'ee par S.-S.~Chern et J.K.~Moser \cite{CM74} pour
classifier les hypersurfaces strictement pseudoconvexes de $\mathbb
C^n$ au moyen d'invariants diff\'erentiels. La complexit\'e des
calculs formels 
ne permet apparemment pas de d\'egager une obstruction
concr\`ete de non alg\'ebrisabilit\'e locale, partag\'ee par tous les
\'el\'ements d'une famille de sous-vari\'et\'es de dimension
quelconque dans $\mathbb C^n$.

On pr\'esente dans cette Note des conditions n\'ecessaires
d'alg\'ebrisabilit\'e locale de sous-vari\'et\'es analytiques
r\'eelles g\'en\'eriques de $\mathbb C^n$. L'argument cl\'e consiste
\`a observer que le (pseudo)groupe des biholomorphismes locaux de
$\mathbb C^n$ stabilisant une sous-vari\'et\'e g\'en\'erique
alg\'ebrique r\'eelle est, sous certaines hypoth\`eses de non
d\'eg\'en\'erescence, un groupe de Lie local pour lequel la loi de
multiplication est alg\'ebrique.

\section{\'Enonc\'e des r\'esultats}

Soit $M$ une sous-vari\'et\'e analytique r\'eelle locale de
codimension $d$ dans $\mathbb C^n$, passant par l'origine, d\'efinie
par les \'equations $r_k(t,\bar{t})=0$ pour $k=1,\dots,d$, o\`u
$t=(t_1,\dots,t_n) \in \mathbb C^n$ et o\`u les $r_k(t,\bar{t}) \in
\mathbb C\{t,\bar{t}\}$ satisfont $r_k(0,0)=0$, $r_k(t,\bar{t}) \equiv
\overline{r_k(t,\bar{t})}$ et $dr_1\wedge \cdots \wedge
dr_d(0) \neq 0$. On suppose $M$ {\sl g\'en\'erique}, c'est-\`a-dire
que $T_0M+iT_0M=T_0\mathbb C^n$. Dans ce cas, $T^{(0,1)}M:=T_0M \cap
T^{(0,1)}\mathbb C^n$ est un fibr\'e complexe de rang $(n-d)$. Soit
$\overline{L}_1,\dots, \overline{L}_{n-d}$ une base de
$T^{(0,1)}M$. On dit que $M$ est {\sl minimale en $0$} si l'orbite (au
sens de Sussmann) de l'origine sous l'action des champs ${\rm
Re}\,\overline{L}_1, {\rm Im}\,\overline{L}_1, \dots, {\rm
Re}\,\overline{L}_{n-d}, {\rm Im}\,\overline{L}_{n-d}$ contient un
voisinage de l'origine dans $M$.  On dit que $M$ est {\sl
holomorphiquement non d\'eg\'en\'er\'ee} ({\it cf}. \cite{Sta96}) s'il
n'existe pas de champ de vecteurs non nul de type $(1,0)$, \`a
coefficients holomorphes, tangent \`a $M$. Dans ce cas, d'apr\`es
\cite{BER99}, il existe un sous-ensemble analytique r\'eel strict $V$
de $M$ tel que, pour tout point $p$ de $M \backslash V$, on a ${\rm
Vect}\, \{\overline{L}^\beta\, \nabla_t (r_k)(p,\bar{p}): \,
\beta \in \mathbb N^{n-d}, \, k=1,\dots,d\}=\C^n$, o\`u l'on note $\nabla_t
(r_k)(t,\bar t)$ le gradient holomorphe de $r_k$ et
$\overline{L}^\beta:=(\overline{L}_1)^{\beta_1}\dots
(\overline{L}_{n-d})^{\beta_{n-d}}$.  On dit alors que $M$ est {\sl finiment
non d\'eg\'en\'er\'ee en $p$}.

On \'etudie la structure du (pseudo)groupe des biholomorphismes locaux
de $\mathbb C^n$ stabilisant $M$. Ce groupe est de dimension infinie
lorsque $M$ n'est nulle part minimale ou holomorphiquement
d\'eg\'en\'er\'ee ({\it cf}. \cite{Sta96}). La notion de {\sl groupe
de Lie local} (de dimension finie) est pr\'esent\'ee par exemple dans
le chapitre 8 de \cite{Sto00}. On appelle {\sl groupe de Lie local
alg\'ebrique} un groupe de Lie local pour lequel la loi de
multiplication ({\it cf}.~\cite{Sto00}, pp.~176--177) est alg\'ebrique
(au sens de Nash).

\begin{theoreme}
\label{theoreme1}
Soit $M$ une sous-vari\'et\'e g\'en\'erique alg\'ebrique r\'eelle de
$\mathbb C^n$, passant par l'origine, minimale et finiment non
d\'eg\'en\'er\'ee en $0$. Il existe $\varepsilon_0 >0$ tel que pour tout $0 < \varepsilon < \varepsilon_0$, les trois conditions suivantes sont satisfaites~:

\noindent {\bf (a)} Le (pseudo)groupe $G_M$ des biholomorphismes locaux
 de $\mathbb C^n$
d\'efinis sur le polydisque 
$\Delta_n(0,\varepsilon)$ et stabilisant $M$ est un groupe de Lie local
alg\'ebrique r\'eel de dimension finie $m$ ne
d\'ependant que de la g\'eom\'etrie locale de $M$ au voisinage de
l'origine.

\noindent {\bf (b)} Il existe $0<\varepsilon'<\varepsilon$ et une application
{\rm alg\'ebrique} $H_M$, constructible algorithmiquement \`a partir 
des \'equations d\'efinissantes de $M$, d\'efinie au voisinage de l'origine 
dans $\mathbb
C^n \times \mathbb R^m$, \`a valeurs dans $\mathbb C^n$, v\'erifiant
$H_M(t,0) \equiv t$, et telle que toute application holomorphe $h
:\Delta_n(0,\varepsilon')\rightarrow~\Delta_n(0,\varepsilon)$ suffisamment
proche de l'identit\'e, stabilisant $M$, v\'erifie $h = H_M(.,e_h)$ pour un unique $e_h
\in \mathbb R^m$.

\noindent {\bf (c)} L'application $(t,e) \mapsto H_M(t,e)$ d\'efinit 
un groupe de Lie local alg\'ebrique de biholomorphismes alg\'ebriques stabilisant $M$.
\end{theoreme}

La d\'emonstration du Th\'eor\`eme \ref{theoreme1} fait intervenir les jets en $0$ de biholomorphismes locaux, voir \cite{GM02-2}.

On consid\`ere maintenant la famille $\mathcal T_n^d$ des
 sous-vari\'et\'es analytiques r\'eelles g\'en\'eriques de $\mathbb
 C^n$, d\'efinies au voisinage de l'origine, de codimension $d$,
 minimales et finiment non d\'eg\'en\'er\'ees en $0$, dont le groupe
 $G_M$ est commutatif, de dimension $n$. Il existe alors, pour tout
 \'el\'ement $M$ de $\mathcal T_n^d$, un syst\`eme local de
 coordonn\'ees holomorphes centr\'e \`a l'origine $(z,w)=(x+iy,u+iv)
 \in \mathbb C^{n-d} \times \mathbb C^d$ et $d$ fonctions analytiques
 r\'eelles $\varphi_1,\dots,\varphi_d$ s'annulant en $0$, tels que $M$
 soit repr\'esent\'e au voisinage de $0$ par les \'equations
 $v_1=\varphi_1(y),\dots,v_d=\varphi_d(y)$. Dans ce cas, $M$ est
 finiment non d\'eg\'en\'er\'e en $0$ si et seulement s'il existe
 $(n-d)$ multiindices $\beta^1,\dots,\beta^{n-d} \in \mathbb N^{n-d}$ de longueur
 strictement positive et des entiers $1 \leq k_1,\dots,k_{n-d}\leq d$
 tels que l'application r\'eelle \def\theequation{1.1}\begin{equation}
 \psi(y):=\left( {\partial^{\vert \beta^1\vert}\varphi_{k_1}(y) \over
 \partial y^{\beta^1}},\ldots, {\partial^{\vert
 \beta^{n-d}\vert}\varphi_{k_{n-d}}(y) \over \partial y^{\beta^{n-d}}}
 \right)
\end{equation}
est de rang $(n-d)$ \`a l'origine dans $\mathbb R^{n-d}$ ({\it
cf}. \cite{GM02-2} Lemme 3.2). On note $y' \mapsto
(\psi'_1(y'),\dots,\psi'_{n-d}(y'))$ l'application r\'eciproque de
$\psi$.  Le second r\'esultat de cette Note est le th\'eor\`eme
suivant~:

\begin{theoreme}
\label{theoreme2}
Soit $M$ appartenant \`a $\mathcal
T_n^d$, d\'efinie par les \'equations
$v_1=\varphi_1(y),\dots,v_d=\varphi_d(y)$, minimale et
finiment non d\'eg\'en\'er\'ee en $0$. Si $M$ est localement
alg\'ebrisable \`a l'origine, les d\'eriv\'ees partielles
$\partial \psi'_j / \partial y'_l$ sont
alg\'ebriques r\'eelles pour $1\leq j,l\leq n-d$ .
\end{theoreme}


Cette condition d'alg\'ebrisabilit\'e locale s'exprime tr\`es
simplement lorsque $M$ est une hypersurface Levi non d\'eg\'en\'er\'ee
de $\mathbb C^n$, d\'efinie par l'\'equation
$v=\varphi(y)$~: il est \'equivalent de dire que les fonctions
$\partial \psi'_j/\partial y'_l$ sont alg\'ebriques r\'eelles pour $1 \leq j,l\leq n-1$
 et que les d\'eriv\'ees
partielles secondes $\partial^2\varphi/\partial y_j\partial y_l$ sont des fonctions
alg\'ebriques des d\'eriv\'ees partielles premi\`eres $\partial \varphi / \partial y_1,
\dots,\partial \varphi / \partial y_{n-1}$ pour $1\leq j,l\leq n-1$.
Le Th\'eor\`eme \ref{theoreme2} montre qu'un \'el\'ement de $\mathcal
T^d_n$ pour lequel une d\'eriv\'ee $\partial \psi'_j/ \partial y'_l$ n'est pas
alg\'ebrique r\'eelle (ce qui est g\'en\'eriquement le cas au sens de
Baire) n'est pas localement alg\'ebrisable \`a l'origine. Le Corollaire
\ref{corollaire1} exhibe une famille explicite d'exemples de
sous-vari\'et\'es analytiques r\'eelles de $\mathbb C^n$, qui ne sont
pas localement alg\'ebrisables \`a l'origine~:

\begin{corollaire} {\rm \cite{GM02-2}}
\label{corollaire1}\ 
{\bf (a)} Soient $\chi_1,\dots,\chi_{n-1}$ des fonctions analytiques
r\'eelles arbitraires, d\'efinies au voisinage de $0 \in \mathbb
R^{n-1}$. L'hypersurface $M$ de $\mathbb C^n$ d'\'equation
$v=\sum_{k=1}^{n-1}[y_k^2+y_k^6+y_k^9y_1\cdots
y_{k-1}+y_k^{n+8}\chi_k(y_1,\dots,y_{n-1})]$ appartient \`a la famille
$\mathcal T^1_n$. De plus, pour un choix g\'en\'erique (au sens de
Baire) de $\chi_1,\dots,\chi_{n-1}$, $M$ n'est pas localement
alg\'ebrisable \`a l'origine.

{\bf (b)} Les hypersurfaces de $\mathbb C^2$ d'\'equations $
v=\sin(y^2),\ v=\sinh(y^2)$ et $v=\exp(\exp(y)-1)-1$ appartiennent \`a
$\mathcal T^1_2$ et ne sont pas localement alg\'ebrisables \`a l'origine.
\end{corollaire}
Pour v\'erifier l'appartenance de $M$ \`a la famille $\mathcal T_n^1$,
on utilise la th\'eorie analytique des sym\'etries des \'equations aux
d\'eriv\'ees partielles ({\it cf}. \cite{Ol86} Chapitre 2), appliqu\'ee
\`a la g\'eom\'etrie CR dans \cite{Su99, Su01, GM02-1}.

Dans le m\^eme \'etat d'esprit, le th\'eor\`eme suivant donne un
\'eclairage diff\'erent du r\'esultat de \cite{HJY01}~:

\begin{theoreme} {\rm \cite{GM02-2}}
\label{theoreme3}\ 
Soit $M : v=\varphi(z\bar{z})$ une hypersurface analytique r\'eelle,
Levi non d\'eg\'en\'er\'ee dans $\mathbb C^2$, passant par l'origine,
dont l'alg\`ebre de Lie des automorphismes infinit\'esimaux est
engendr\'ee par $\partial / \partial w$ et $iz\partial / \partial z$. Si $M$ est localement
alg\'ebrisable \`a l'origine, la d\'eriv\'ee  de $\varphi$ est 
alg\'ebrique. Ainsi, les hypersurfaces de $\mathbb
C^2$ d'\'equations $v={\rm e}^{z\bar{z}}-1$, $v=\sin(z\bar{z})$ et $
v=\sinh(z\bar{z})$ ne sont pas localement alg\'ebrisables \`a
l'origine.
\end{theoreme}

\section{R\'esum\'e de la d\'emonstration du Th\'eor\`eme \ref{theoreme2}}

Soit $M \in \mathcal T^d_n$ et soit $\Phi$ un biholomorphisme local,
v\'erifiant $\Phi(0)=0$, tel que $M':=\Phi(M)$ est alg\'ebrique,
minimale et finiment non d\'eg\'en\'er\'ee en $0$.  L'alg\`ebre de Lie
r\'eelle des automorphismes infinit\'esimaux de $M$ \'etant
engendr\'ee par les champs $\partial/\partial t_1, \dots, \partial/\partial t_n$,
l'alg\`ebre de Lie des automorphismes infinit\'esimaux de $M'$ est
commutative, de dimension $n$, engendr\'ee par les champs
$X'_1=\Phi_*(\partial/\partial t_1), \dots,
X'_n=\Phi_*(\partial/\partial t_n)$. On commence par redresser
alg\'ebriquement les feuilletages holomorphes induits par
$X'_1,\dots,X'_n$, en utilisant l'alg\'ebricit\'e de l'application
$H_{M'}$ et la commutativit\'e du groupe $G_{M'}$ donn\'ees dans le
Th\'eor\`eme \ref{theoreme1}. La d\'emonstration est standard, voir \cite{GM02-2}.

\begin{proposition}
\label{proposition1}
 Il existe un biholomorphisme alg\'ebrique complexe $\Psi$, d\'efini au
voisinage de l'origine dans $\mathbb C^n$, et $n$ fonctions
alg\'ebriques complexes $c_1,\dots,c_n$, v\'erifiant $c_j(0) = 1$,
tels que l'alg\`ebre de Lie des automorphismes infinit\'esimaux de
$M'':=\Psi(M')$ est engendr\'ee par les champs
$X''_j:=\Psi_*(X'_j)=c_j(t_j)\partial/\partial t_j$, $1 \leq j \leq n$.
\end{proposition}

D'apr\`es \cite{CM74}, il existe des coordonn\'ees $t=(z_1,\dots,z_{n-d},w_1,\dots,w_d)$ dans lesquelles on peut repr\'esenter $M''$ par les \'equations
$w_k=\overline{\Theta}_k(z,\bar{z},\bar{w})$, $k=1,\dots,d$,
o\`u $\overline{\Theta}_k$ est une fonction alg\'ebrique complexe. On note maintenant
, $(\Psi \circ
\Phi)^{-1}=(f_1,\dots,f_{n-d},g_1,\dots,g_d)$ et $(c_1,\dots,c_n) =
(a_1,\dots,a_{n-d},b_1,\dots,b_d)$. L'application $\Psi \circ \Phi$
v\'erifiant $(\Psi \circ \Phi)_*(\partial/\partial t_j) = c_j(t_j)
\partial/\partial t_j$, on a les identit\'es
$a_j(z_j)(\partial f_j/\partial z_j)(z_j) \equiv b_k(w_k)(\partial g_k/\partial w_k)(w_k) \equiv 1$,
ce qui montre que $\partial f_j/\partial z_j$ et $\partial g_k/\partial w_k$ 
sont alg\'ebriques pour $1 \leq j \leq n-d$
et $1 \leq k \leq d$.  La
condition $(\Psi \circ \Phi)(M) = M''$ donne, pour
$k=1,\dots,d$, les identit\'es suivantes dans $\mathbb C\{z,\bar{z},\bar{w}\}$~:\def\theequation{1.2}\begin{equation}
\frac{g_k(\overline{\Theta}_k(z,\bar{z},\bar{w}))-
\bar{g}_k(\bar{w}_k)}{2i}  \equiv  
\varphi_k\left(\frac{f_1(z_1)-\bar{f}_1(\bar{z}_1)}{2i},
\dots,\frac{f_{n-d}(z_{n-d})-\bar{f}_{n-d}(\bar{z}_{n-d})}{2i}\right).
\end{equation}
En diff\'erentiant (1.2) par rapport \`a $z_j$ pour
$j=1,\dots,n-d$, on obtient~: 
\def\theequation{1.3}\begin{equation}
\frac{a_j(z_j)(\partial \overline{\Theta}_k/\partial z_j)
(z,\bar{z},\bar{w})}{b_k(\overline{\Theta}_k(z,\bar{z},\bar{w}))}
\equiv \frac{\partial \varphi_k}{\partial
y_j}\left(\frac{f_1(z_1)-\bar{f}_1(\bar{z}_1)}{2i},
\dots,\frac{f_{n-d}(z_{n-d})-\bar{f}_{n-d}(\bar{z}_{n-d})}{2i}\right).
\end{equation}
Le membre de gauche de l'\'equation (1.3) est une fonction
alg\'ebrique complexe (l'alg\'ebricit\'e est stable par diff\'erentiation, {\it cf}. \cite{BCR98}), ind\'ependante de la variable $w$. Pour $1 \leq l \leq n-d$, on applique l'op\'erateur $\partial^{|\beta^l|}/ \partial z^{\beta^l}$ \`a l'identit\'e (1.2) avec $k=k_l$. 
Il existe des fonctions
alg\'ebriques r\'eelles $A_{k_l,\beta^l}$ v\'erifiant~:
\def\theequation{1.4}\begin{equation}
A_{k_l,\beta^l}(z,\bar{z}) \equiv
\frac{\partial^{|\beta^l|} \varphi_{k_l}}{\partial
y^{\beta^l}}\left(\frac{f_1(z_1)-\bar{f}_1(\bar{z}_1)}{2i},
\dots,\frac{f_{n-d}(z_{n-d})-\bar{f}_{n-d}(\bar{z}_{n-d})}{2i}\right).
\end{equation}
L'application $\psi$ donn\'ee par l'\'equation (1.1) \'etant de rang
$n-d$, on obtient pour $j=1,\dots,n-d$ l'identit\'e
$f_{j}(z_{j})-\bar{f}_{j}(\bar{z}_{j}) \equiv
2i\,\psi'_j(A_{k_1,\beta^1}(z,\bar{z}),\dots,
A_{k_{n-d},\beta^{n-d}}(z,\bar{z}))$. La diff\'erentiation de cette
identit\'e par rapport \`a $z_j$ et \`a $z_m$ pour $m\neq j$ fournit le syst\`eme lin\'eaire~:
\def\theequation{1.5}\begin{equation}
\left\{
\aligned
{1\over 2i\, a_j(z_j)}\equiv 
& \ 
\sum_{l=1}^{n-d}\, 
{\partial \psi_j'\over\partial y_l'}
({A}_{k_1,\beta^1}(z,\bar z),\dots,
{A}_{k_{n-d},\beta^{n-d}}(z,\bar z))\,
{\partial {A}_{k_l,\beta^l}
\over \partial z_j}
(z,\bar z),\\
0 \equiv
& \ 
\sum_{l=1}^{n-d}\, {\partial\psi_j'\over\partial y_l'}
({A}_{k_1,\beta^1}(z,\bar z),\dots,{A}_{k_{n-d},\beta^{n-d}}(z,\bar z))\,
{\partial {A}_{k_l,\beta^l}\over\partial z_{m}}(z,\bar z), \ \ \ \ \
m\neq j.
\endaligned\right.
\end{equation}
Puisque $a_j(0) =1$ pour tout $1 \leq j \leq n-d$, la matrice $((\partial A_{k_l,\beta^l}/\partial z_j) (0,0))_{1 \leq j,l \leq
n-d}$ est inversible et il existe donc des fonctions
alg\'ebriques complexes $B_{j,l}(z,\bar{z})$, d\'efinies pour $1 \leq
j,l\leq n-d$, telles que $(\partial \psi'_j/\partial y_l')
(A(z,\bar{z})):= (\partial \psi'_j/\partial y_l')
(A_{k_1,\beta^1}(z,\bar{z}),\dots,A_{k_{n-d},\beta^{n-d}}(z,\bar{z}))
\equiv B_{j,l}(z,\bar{z})$.  Enfin, la fonction $A_{k_l,\beta^l}$
\'etant r\'eelle et la matrice $((\partial A_{k_l,\beta^l}/\partial z_j) (0,0))_{1 \leq j,l
\leq n-d}$ \'etant inversible, le jacobien \`a l'origine de
l'application $y \mapsto A(iy,-iy)=:y'$ est non nul. Il existe donc
une application alg\'ebrique r\'eelle $C$ telle
que $(\partial \psi'_{j}/\partial y'_l)(y') \equiv
B_{j,l}(iC(y'),-iC(y'))$, ce qui prouve l'alg\'ebricit\'e de $\partial \psi'_{j}/\partial y'_l$ pour $1 \leq j,l\leq n-d$. \qed

\label{lastpage}


\begin{thebibliography}{99}
\renewcommand{\partopsep}{0pt}
\renewcommand{\itemsep}{-1pt}
\renewcommand{\parsep}{0pt}

\bibitem{BER99}
Baouendi, M.S., Ebenfelt, P., Rothschild, L.P., Rational dependence of
smooth and analytic CR mappings on their jets, Math. Ann.  315 (1999),
205--249.

\bibitem{BCR98}
Bochnak, J., Coste, M., Roy, M.F., Real algebraic geometry, 
Springer-Verlag, Berlin, 1998. x+430~pp.

\bibitem{Ca32} Cartan, \'E., Sur la g\'eom\'etrie pseudo-conforme des
hypersurfaces de l'espace de deux variables complexes, I, Annali di
Mat. 11 (1932), 17--90, II, Annali Sc. Norm. Sup. Pisa 1 (1932),
333--354.

\bibitem{CM74}
Chern, S.S., {\sc Moser}, J.K., Real hypersurfaces in complex manifolds,
Acta Math. 133 (1974), no.~2, 219--271.

\bibitem{GM02-1} 
Gaussier, H., Merker, J., A new example of uniformly Levi degenerate
hypersurface in $\mathbb C^3$, Ark. Mat. (to appear).

\bibitem{GM02-2} Gaussier, H., Merker, J., Nonalgebraizable
real analytic tubes in $\mathbb C^n$, Pr\'epublication LATP (2002), 
no.~02--17, 36~pp. 

\bibitem{HJY01} Huang, X., Ji, S., Yau, S.S., An example of a real
analytic strongly pseudoconvex hypersurface which is not
holomorphically equivalent to any algebraic hypersurface,
Ark. Mat. 39 (2001), no.~1, 75--93.

\bibitem{Ol86} Olver, P.J., Applications of Lie groups to
differential equations, Springer Verlag, Heidelberg, 1986.

\bibitem{Sta96} Stanton, N., Infinitesimal CR automorphisms of real hypersurfaces, Amer. J. Math.  118  (1996),  no.~1, 209--233.

\bibitem{Sto00} Stormark, O., Lie's structural approach to PDE systems, 
Cambridge University Press, Cambridge, 2000.

\bibitem{Su99} Sukhov, A., Segre varieties and Lie
symmetries, Math. Z. 231 (2001), no.~3, 483--492.

\bibitem{Su01} Sukhov, A., On maps of CR manifolds and transformations
of differential equations, C. R. Acad. Sci. Paris, S\'erie~I 333
(2001), 545--550.




\end{thebibliography}
\end{document}